\documentclass[12pt]{amsart}
\usepackage{amssymb,amsfonts,amsmath}

\begin{document}

\title[Amalgamated free product over Cartan subalgebra, II]
{Amalgamated free product over Cartan Subalgebra, II \\
Supplementary Results \& Examples}

\author[Y. Ueda]{Yoshimichi UEDA}
\address{Department of Mathematics, Graduate School of Science, Hiroshima University, 
Higashi-Hiroshima, 739-8526, Japan}
\address{
{\it Current Address}: Faculty of Mathematics, Kyushu University, Fukuoka, 810-8560, 
Japan}
\email{e-mail: ueda@math.kyushu-u.ac.jp} 
\maketitle

\newtheorem{thm}{Theorem}[section]
\newtheorem{cor}[thm]{Corollary}
\newtheorem{lemma}[thm]{Lemma}
\newtheorem{prop}[thm]{Proposition}
\newtheorem{fact}[thm]{Fact}
\newtheorem{remark}[thm]{Remark}
\newtheorem{remarks}[thm]{Remarks}

\theoremstyle{definition}
\newtheorem{defn}{Definition}[section]
\baselineskip=15pt

\section{Introduction}

Let $A \supseteq D \subseteq B$ be two von Neumann algebras 
together with a common Cartan subalgebra. 
Then the amalgamated free product $M = A *_D B$ 
with respect to the unique conditional expectations 
from $A$, $B$ onto $D$ can be considered. 
In our previous paper [U1], the questions of 
its factoriality and type classification were discussed 
in detail, which will be reviewed in \S4. The main purpose of 
the paper is to give further  supplementary results obtained after 
the completion of the previous paper together with discussing some
examples. 

\medskip
The author would like to express his sincere gratitude to 
the organizers Bruce Blackadar \& Hideki Kosaki for 
inviting him to the US-Japan seminar 1999 held at Fukuoka, 
Japan and for giving this opportunity to present this paper. 

\section{Amalgamated Free Products of von Neumann algebras}

Let $A \supseteq D \subseteq B$ be $\sigma$-finite 
von Neumann algebras, and let 
$E_D^A : A \rightarrow D$, 
$E_D^B : B \rightarrow D$ 
be faithful normal conditional expectations. 
Then one can consider the amalgamated free product of 
$A$ and $B$ over $D$ with respect to the conditional 
expectations $E_D^A$, $E_D^B$:  
$$
(M, E_D^M) = (A, E_D^A) *_D (B, E_D^B),  
$$
which is defined as a pair of a von Neumann algebra $M$ 
into which the triple $A \supseteq D \subseteq B$ is embedded 
and a faithful normal conditional expectation 
$E_D^M : M \rightarrow D$, and characterized by 
the following three conditions: 
\begin{itemize}
\item $M$ is generated by the subalgebras $A$, $B$; 
\item $E_D^M|_A = E_D^A$, $E_D^M|_B = E_D^B$; 
\item $A$, $B$ are free with amalgamation over $D$ 
                 in the $D$-probability space 
                 $(M \supseteq D, E_D^M)$, see [VDN], i.e.,   
                 $$
                 E_D^M(
                 \{\text{alternating words in 
                 $A^{\circ}, B^{\circ}$}\}
                 ) = 0, 
                 $$ 
                 where we denote $A^{\circ} = \text{Ker}E_D^A$, 
                 $B^{\circ} = \text{Ker}E_D^B$ as usual. 
\end{itemize} 
For the details, we refer to [P2],[VDN],[U1] 
(see also [BD]). 

\medskip
In analysis on type III factors, modular automorphisms 
are of central importance so that we need 
to compute the modular automorphisms 
$\sigma_t^{\varphi\circ E_D^M}$ ($t \in {\mathbf R}$) 
for a faithful normal state $\varphi$  on $D$. 

\begin{thm}
 {\bf ([U1, Theorem 2.6])} 
We have 
\begin{equation}
\sigma_t^{\varphi\circ E_D^M}|_A = 
\sigma_t^{\varphi\circ E_D^A}, \quad 
\sigma_t^{\varphi\circ E_D^M}|_B = 
\sigma_t^{\varphi\circ E_D^B}. \tag{2.1}
\end{equation}
\end{thm} 

The modular operator $\Delta_{\varphi\circ E_D^M}$ and 
the modular conjugation $J^M$ can be also computed 
explicitly. (See [U1, Appendix I].) 

\section{Amalgamated Free Products over Cartan Subalgebras}

Let $A$ and $B$ be von Neumann algebras with separable 
preduals, and we suppose that they have a common Cartan 
subalgebra $D$ or equivalently that there is a common 
subalgebra $D$ satisfying: 
\begin{itemize}
\item $D$ is a MASA in both $A$ and $B$; 
\item there are (automatically unique faithful) normal 
                 conditional expectations 
                 $E_D^A : A \rightarrow D$, 
                 $E_D^B : B \rightarrow D$; 
\item the normalizers ${\mathcal N}_A(D)$, ${\mathcal N}_B(D)$ 
                 generate the whole $A$, $B$, respectively. 
\end{itemize}
(See [FM].) Let 
$$
(M, E_D^M) = (A, E_D^A) *_D (B, E_D^B), 
$$
and we will write $M = A *_D B$ since there is no other 
choice of normal conditional expectations from $A$, $B$ 
onto $D$. 

\medskip
The triple $A \supseteq D \subseteq B$ produces two 
countable non-singular Borel equivalence relations 
${\mathcal R}_A$, ${\mathcal R}_B$ over a common standard Borel 
probability space $(X,\mu)$ in such a way that 
\begin{equation}
A = W^*_{\sigma_A}({\mathcal R}_A), \quad 
B = W^*_{\sigma_B}({\mathcal R}_B), \quad 
D = L^{\infty}(X,\mu), \tag{3.1}
\end{equation}
where $W^*_{\sigma_A}({\mathcal R}_A)$, 
$W^*_{\sigma_B}({\mathcal R}_B)$ denote 
the von Neumann algebras constructed from 
${\mathcal R}_A$, 
${\mathcal R}_B$ 
together with relevant 2-cocycles 
$\sigma_A$, $\sigma_B$, respectively, 
by the Feldman-Moore construction ([FM]). 
After fixing a point  realization 
$D = L^{\infty}(X,\mu)$, 
such a pair (${\mathcal R}_A$, ${\mathcal R}_B$) is uniquely 
determined up to null set. (This fact will be discussed 
in [U2] in detail.) Therefore, the countable non-singular 
Borel equivalence relation 
$$
{\mathcal R}_M := {\mathcal R}_A \vee {\mathcal R}_B \ 
(\subseteq X\times X)
$$
is a canonical object attached to the triple 
$A \supseteq D \subseteq B$, and we call 
this equivalence relation the 
{\it canonical equivalence relation} associated with 
the amalgamated free product $M = A *_D B$ or 
the triple $A \supseteq D \subseteq B$. 

\section{Factoriality \& Type Classification ([U1])}

In this section, we review our previous paper [U1].  

\medskip
Let $A \supseteq D \subseteq B$ be as in \S3, i.e., 
two von Neumann algebras (with separable preduals) and 
a common Cartan subalgebra. 
We here discuss the amalgamated free product 
$M = A *_D B$. 
The first problem is its factoriality, namely, 
to find a suitable sufficient condition for 
the amalgamated free product $M$ to be a factor. 
A satisfactory answer to the problem was given 
in our previous paper. 

\begin{thm}
{\bf ([U1, Theorem 4.3])} 
If either $A$ or $B$ is a factor of non-type {\rm I}, 
then the amalgamated free product $M = A *_D B$ over 
a common Cartan subalgebra $D$ becomes a factor. 
More precisely, if $A$ {\rm (}or $B${\rm )} is a factor 
of non-type {\rm I}, then there is a faithful normal 
state $\varphi$ on $D$ such that 
\begin{equation}
(A_{\varphi\circ E_D^A})' \cap M \subseteq A \quad 
\left(\text{resp.}\ 
(B_{\varphi\circ E_D^B})' \cap M \subseteq B
\right). \tag{4.1}
\end{equation}
Furthermore, if $A$ {\rm (}or $B${\rm )} is further 
assumed to be of type {\rm III$_{\lambda}$} 
$(0 < \lambda \leq 1)$, then the state $\varphi$ can be 
chosen in such a way that 
\begin{equation}
(A_{\varphi\circ E_D^A})' \cap A = {\mathbf C}1 \quad 
\left(\text{resp.}\ 
(B_{\varphi\circ E_D^B})' \cap B = {\mathbf C}1 
\right). \tag{4.2}
\end{equation}
\end{thm} 

This result can be generalized further. 
Such a generalization will be discussed later, 
see Remark 4.8 (2),(3). 

\medskip
The second problem seems 
Murray-von Neumann-Connes' type classification 
of the amalgamated free product $M$. 
In this direction, we obtained 
the following corollaries of Theorem 4.1: 

\begin{cor}
{\bf [U1, Corollary 4.5])} 
Suppose that both $A$ and $B$ are factors 
of non-type {\rm I}. If $M = A *_D B$ is of 
type {\rm III$_0$}, then both $A$ and $B$ must 
also be of type {\rm III$_0$}. 
\end{cor} 

\begin{cor}
{\bf ([U1, p.377])} 
Suppose that both $A$ and $B$ are factors 
of non-type {\rm I}. 

\noindent
{\rm (1)} If either $A$ or $B$ of type {\rm II$_1$} 
and if $M = A *_D B$ is semi-finite {\rm (}i.e., has 
a faithful semi-finite normal trace{\rm )}, then 
$M$ must be of type {\rm II$_1$}. 

\noindent
{\rm (2)} If either $A$ or $B$ is of type 
{\rm III$_{\lambda}$} $(0 < \lambda < 1)$, then 
$M = A *_D B$ must be 
of type {\rm III$_{\lambda^{1/n}}$} or 
of type {\rm III$_1$}. 

\noindent
{\rm (3)} If either $A$ or $B$ is 
of type {\rm III$_1$}, then $M = A *_D B$ must be 
of type {\rm III$_1$}. 

\noindent
{\rm (4)} If $A$ is of type {\rm III$_{\lambda}$} and 
$B$ of type {\rm III$_{\mu}$} with 
$\displaystyle{\frac{\log\lambda}{\log\mu} 
\not\in {\mathbb Q}}$, then $M = A *_D B$ must be 
of type {\rm III$_1$}. 

\end{cor}   

\bigskip
For a while, we assume that $A \supseteq D \subseteq B$ 
is a general triple of $\sigma$-finite von Neumann algebras 
together with faithful normal conditional expectations 
$E_D^A : A \rightarrow D$, 
$E_D^B : B \rightarrow D$, and let  
$$
(M, E_D^M) = (A, E_D^A) *_D (B, E_D^B)
$$
be the amalgamated free product. 
Choose and fix a faithful normal state $\varphi$ on $D$, 
and we set: 
$$
\widetilde{A} := 
A \rtimes_{\sigma^{\varphi\circ E_D^A}} {\mathbf R} 
\supseteq 
\widetilde{D} := 
D \rtimes_{\sigma^{\varphi}} {\mathbf R} 
\subseteq 
\widetilde{B} := 
B \rtimes_{\sigma^{\varphi\circ E_D^B}} {\mathbf R}, 
$$
and 
$$
\widetilde{M} : = 
M \rtimes_{\sigma^{\varphi\circ E_D^M}} {\mathbf R}. 
$$
Then there are faithful normal conditional expectations 
$$
\begin{aligned}
\widehat{E_D^A} : \widetilde{A} \rightarrow 
\widetilde{D}\ ; \quad \quad 
&\widehat{E_D^A}\left(\int_{-\infty}^{\infty} 
a(t)\lambda(t) dt\right) 
:= \int_{-\infty}^{\infty} E_D^A(a(t)) \lambda(t) dt, 
\\
\widehat{E_D^B} : \widetilde{B} \rightarrow 
\widetilde{D}\ ; \quad \quad 
&\widehat{E_D^B}\left(\int_{-\infty}^{\infty} 
b(t)\lambda(t) dt\right) 
:= \int_{-\infty}^{\infty} E_D^B(b(t)) \lambda(t) dt, 
\\ 
\widehat{E_D^M} : \widetilde{M} \rightarrow 
\widetilde{D}\ ; \quad \quad 
&\widehat{E_D^M}\left(\int_{-\infty}^{\infty} 
m(t)\lambda(t) dt\right) 
:= \int_{-\infty}^{\infty} E_D^M(m(t)) \lambda(t) dt. 
\end{aligned}
$$

\begin{thm}
{\bf ([U1, Theorem 5.1])} 
In the current general setting, we have
\begin{equation}
\left(\widetilde{M}, \widehat{E_D^M}\right) \cong 
\left(\widetilde{A}, \widehat{E_D^A}\right) 
*_{\widetilde{D}} 
\left(\widetilde{B}, \widehat{E_D^B}\right). \tag{4.3} 
\end{equation}
Moreover, the dual action $\theta_t^M$ 
$(t \in {\mathbf R})$ associated with $M$ is determined 
by those 
$\theta_t^A$, 
$\theta_t^B$ 
associated with $A$, $B$, respectively, 
in such a way that  
\begin{equation}
\theta_t^M|_{\widetilde{A}} = \theta_t^A, \quad 
\theta_t^M|_{\widetilde{B}} = \theta_t^B. \tag{4.4} 
\end{equation}
\end{thm} 

\bigskip
Let us return to the original setting, namely, 
the triple $A \supseteq D \subseteq B$ consists of 
two von Neumann algebras with separable preduals
and a common Cartan subalgebra. By Theorem 4.4, 
we have 
\begin{equation}
\widetilde{M} \left(= \widetilde{A *_D B}\right) \cong 
\widetilde{A} *_{\widetilde{D}} \widetilde{B}. \tag{4.5}
\end{equation}
Since $\widetilde{D}$ is also a common Cartan subalgebra 
in both $\widetilde{A}$ and $\widetilde{B}$, we write 
$\widetilde{A} *_{\widetilde{D}} \widetilde{B}$ as the 
amalgamated free product von Neumann algebra of 
$\widetilde{A}$ and 
$\widetilde{B}$ over 
$\widetilde{D}$ with respect to 
the conditional expectations $\widehat{E_D^A}$, $\widehat{E_D^B}$ 
since no confusion is possible.   

\medskip
We further suppose that both $A$ and $B$ are factors 
of non-type I in what follows. Theorem 4.4 (or (4.5)) 
together with the proof of Theorem 4.2 implies 

\begin{thm}
{\bf ([U1, Theorem 5.4])} 
In the current setting, 
we have
\begin{equation}
{\mathcal Z}(\widetilde{M}) = 
{\mathcal Z}(\widetilde{A}) \cap {\mathcal Z}(\widetilde{B}) 
\subseteq 
\widetilde{D}. \tag{4.6}
\end{equation}
\end{thm}

Let $(X_A, F^A_t)$, $(X_B, F^B_t)$ be 
the flows of weights ([CT]) of $A$, $B$, respectively. 
Fix a point realization $D =  L^{\infty}(X,\mu)$, 
and set $X_D := X \times {\mathbf R}$ equipped with 
the usual product measure 
$d\mu\otimes e^{-t}dt$ and $F^D_t(x,s) := (x,s+t)$. 
Then there are two factor maps 
$$
\pi^D_A : (X_D, F_t^D) \rightarrow (X_A, F_t^A), \quad 
\pi^D_B : (X_D, F^D_t) \rightarrow (X_B, F^B_t) 
$$
since $\widetilde{D}$ is a common Cartan subalgebra 
in both $\widetilde{A}$, $\widetilde{B}$. 
Let $(X_M, F^M_t)$ be the flow of weights of $M$. 
Theorem 4.5 says that there are three factor maps
$$
\pi^A_M : (X_A, F_t^A) \rightarrow (X_M, F_t^M), \quad 
\pi^B_M : (X_B, F^D_t) \rightarrow (X_M, F^M_t),  
$$
$$
\pi^D_M : (X_D, F_t^D) \rightarrow (X_M, F_t^M).  
$$

\begin{cor}
{\bf ([U1, Corollary 5.6])} 
The flow $(X_M, F^M_t)$ is determined as 
the unique maximal common factor flow 
of those 
$(X_D, F_t^D)$, $(X_A, F^A_t)$, and $(X_B, F^B_t)$. 
\end{cor} 

This corollary explains all the type classification 
results mentioned before. Indeed, if $A$ (or $B$) is 
a factor of type III$_1$, then the flow of weights 
$(X_A, F^A_t)$ (resp. $(X_B, F^B_t)$) is trivial 
(see [CT],[T2]), i.e., the one-point flow, and hence 
the corollary says that so is the flow of weights 
$(X_M, F^M_t)$, which means that $M$ is of type III$_1$. 
The others can be also explained similarly.  

\begin{cor}
{\bf ([U1, Corollary 5.8])} 
If $A$ and $B$ coincide with each other, i.e., $A = B$, 
then the amalgamated free product 
$M = A *_D B$ $(= A *_D A)$ and $A\ (= B)$ have 
the same flow of weights. In particular, 
any ergodic flow can be realized as 
the flow of weights of a certain amalgamated 
free product. 
\end{cor}

\medskip\noindent
\begin{remarks} {\rm A few remarks are in order. 

\noindent
(1) Corollary 4.6 has the trivial reformulation: 
The flow of weights $(X_M, F^M_t)$ coincides with 
the associated flow ([HOO],[Kr],[FM]) of 
the canonical equivalence relation ${\mathcal R}_M$ 
introduced in \S3.  

\noindent
(2) Based on Theorem 4.4, (4.3) 
(or (4.5)) together with the proof of 
Theorem 4.1, we can show the following: 
Let $A$ and  $B$ be von Neumann algebras 
(with separable preduals) having no type I direct 
summand, and let $D$ be a common Cartan subalgebra.  
Then we have 
\begin{equation}
{\mathcal Z}(\widetilde{A *_D B}) = 
{\mathcal Z}(\widetilde{A})\cap{\mathcal Z}(\widetilde{B})\ 
(\subseteq \widetilde{D}).  \tag{4.7}
\end{equation}
This implies 
\begin{equation}
{\mathcal Z}(A *_D B) = 
{\mathcal Z}(\widetilde{A *_D B})^{\theta^M} = 
{\mathcal Z}(\widetilde{A})^{\theta^A}\cap{\mathcal 
Z}(\widetilde{B})^{\theta^B} 
= {\mathcal Z}(A) \cap {\mathcal Z}(B)\ 
(\subseteq D) \tag{4.8}
\end{equation}
thanks to (4.4) and 
the continuous decomposition theorem [T2].  

\noindent
(3) The (4.8) in the above (2) can be 
reformulated as follows: 
Under the assumption that both $A$ and $B$ have 
no type I direct summand, 
the amalgamated free product $M = A *_D B$ 
is a factor if and only if 
the canonical equivalence relation 
${\mathcal R}_M$ is ergodic.
}    
\end{remarks}

\section{Miscellaneous Results}

\medskip\noindent
{\bf 5.1.}\ 
Let $A \supseteq D \subseteq B$ and 
$E_D^A : A\rightarrow D$, 
$E_D^B : B \rightarrow D$ 
be as in \S2, and let 
$$
(M, E_D^M) = (A, E_D^A) *_D (B, E_D^B) 
$$
be the amalgamated free product. Suppose that $M$ has 
separable predual, or equivalently that so do both $A$ 
and $B$, and further that we have known that 
\begin{equation}
{\mathcal Z}(M) = {\mathcal Z}(A) \cap {\mathcal Z}(B) \subseteq 
{\mathcal Z}(D). \tag{5.1.1} 
\end{equation}
Then we have the following simultaneous direct 
integral decompositions: 
$$
M = \int_{\Omega}^{\oplus} M(\omega)\ d\nu(\omega) \ 
\supseteq\  
A = \int_{\Omega}^{\oplus} A(\omega)\ d\nu(\omega)\  
\supseteq\  
D = \int_{\Omega}^{\oplus} D(\omega)\ d\nu(\omega) 
$$
$$
M = \int_{\Omega}^{\oplus} M(\omega)\ d\nu(\omega) \ 
\supseteq 
B = \int_{\Omega}^{\oplus} B(\omega)\ d\nu(\omega)\   
\supseteq\  
D = \int_{\Omega}^{\oplus} D(\omega)\ d\nu(\omega) 
$$
with ${\mathcal Z}(M) = L^{\infty}(\Omega,\nu)$. Let 
$$
{\mathcal H} = \int_{\Omega}^{\oplus} {\mathcal H}(\omega)\ 
d\nu(\omega)
$$
be the corresponding direct integral decomposition 
of ${\mathcal H} =L^2(M)$. We may and do assume that 
$\omega \mapsto {\mathcal H}(\omega)$ is a constant field 
of the separable infinite dimensional Hilbert space. 
The conditional expectations $E_D^M$, $E_D^A$, $E_D^B$ 
are also decomposed as follows: 
$$
E_D^M = \int_{\Omega}^{\oplus} (E_D^M)_{\omega}\ 
d\nu(\omega), 
$$
$$
E_D^A = \int_{\Omega}^{\oplus} (E_D^A)_{\omega}\ 
d\nu(\omega), \quad 
E_D^B = \int_{\Omega}^{\oplus} (E_D^B)_{\omega}\ 
d\nu(\omega), 
$$ 

\medskip
\begin{thm}
Under the hypothesis (5.1.1), we have
\begin{equation}
(M(\omega), (E_D^M)_{\omega}) \cong 
(A(\omega), (E_D^A)_{\omega}) *_{D(\omega)} 
(B(\omega), (E_D^B)_{\omega}) \tag{5.1.2} 
\end{equation}
for almost every $\omega \in \Omega$. 
\end{thm}

Before going to the proof, we provide a suitable 
(for our purpose) reformulation of freeness. 
Let $(N \supseteq L,\ E : N \rightarrow L)$ be 
a $L$-probability space, i.e., $N \supseteq L$ is 
an inclusion of unital algebras with the same unit and $E : N
\rightarrow  L$ is a conditional expectation 
in the purely algebraic sense.  
Assume that $N_1$, $N_2$ are unital algebras 
containing $L$ in common. We introduce
the operation 
$x \in N \mapsto [x]^{\circ} := x - E(x)$, and 
the freeness (with amalgamation over $L$) of the pair 
$N_1$, $N_2$ can be interpreted as follows: 

\begin{lemma}
The pair $N_1$, $N_2$ are free with amalgamation 
over $L$ if and only if the mapping
$$
\Phi_{(N_1, N_2)} : (x_1, x_2, \cdots, x_n) \in 
\Lambda(N_1, N_2) 
\mapsto E([x_1]^{\circ} [x_2]^{\circ} \cdots 
[x_n]^{\circ}) 
$$
is identically zero. Here, $\Lambda(N_1, N_2)$ is 
the set of those finite alternating sequences 
$(x_1, x_2, \cdots, x_n)$ of elements 
in $N_1 \cup N_2$ 
with $x_i \in N_{j(i)}$,\ 
$j(1) \neq j(2) \neq \cdots \neq j(n)$. 
\end{lemma}

In our case, the operation 
$m \in M \mapsto [m]^{\circ} 
:= m - E_D^M(m) 
\in M^{\circ} := \text{Ker}E_D^M$ 
is normal and linear, 
and can be 
(direct integral) decomposed as follows: 
$$
[\ \cdot\ ]^{\circ} = 
\int_{\Omega}^{\oplus} [\ \cdot\ ]^{\circ}_{\omega}\ 
d\nu(\omega),   
$$
where $[m(\omega)]^{\circ}_{\omega} = 
m(\omega) - (E_D^M)_{\omega}(m(\omega))$, 
a normal linear map, for almost every 
$\omega \in \Omega$. 

\medskip
\begin{proof}(Proof of Theorem 5.1)\ Since both $A$ and $B$ have 
separable preduals, one can choose countable families 
$\{a_k\}_{k \in {\mathbb N}}$, 
$\{b_k\}_{k \in {\mathbb N}}$ 
in such a way that they generate $A$ and $B$, respectively. 
Let 
$$
a_k = \int_{\Omega}^{\oplus} a_k(\omega)\ d\nu(\omega), 
\quad 
b_k = \int_{\Omega}^{\oplus} b_k(\omega)\ d\nu(\omega). 
$$
Then we can choose a co-null Borel subset 
$\Omega_1$ of $\Omega$ in such a way that,  
for every $\omega \in \Omega_1$, 
\begin{itemize}
\item[(a)] $M(\omega)$ is generated by $A(\omega)$ and $B(\omega)$; 
\item[(b)] 
$(E_D^M)_{\omega} : M(\omega) \rightarrow D(\omega)$, 
$(E_D^A)_{\omega} : A(\omega) \rightarrow D(\omega)$, 
$(E_D^B)_{\omega} : B(\omega) \rightarrow D(\omega)$  
are faithful normal conditional expectations and 
$$
(E_D^M)_{\omega}|_{A(\omega)} = 
(E_D^A)_{\omega}, \quad 
(E_D^M)_{\omega}|_{B(\omega)} = 
(E_D^B)_{\omega}; 
$$
\item[(c)] $A(\omega)$ and $B(\omega)$ are generated by 
the $a_k(\omega)$'s and the $b_k(\omega)$'s, 
respectively. 
\end{itemize}
Replacing the $a_k$'s (resp. the $b_k$'s) by all the finite 
products of them and their adjoints from the beginning, 
we may and do assume, instead of (c), that 
\begin{itemize}
\item[(c)'] the linear span of the $a_k(\omega)$'s 
(resp. the $b_k(\omega)$'s) forms 
a $\sigma$-weakly dense $*$-subalgebra of 
$A(\omega)$ (resp. $B(\omega)$). 
\end{itemize}

The restriction of the map $\Phi_{(A,B)}$ to 
the subset $\Lambda(\{a_k\}_{k\in{\mathbb N}}, 
\{b_k\}_{k\in{\mathbb N}})$ is identically zero by 
the freeness of the pair $A$, $B$ (see Lemma 5.2). 
Therefore, there is a co-null Borel subset 
$\Omega_0$ of $\Omega_1$ such that 
$$
\begin{aligned}
\Phi_{(A(\omega),B(\omega))}&
\left((m_1(\omega),m_2(\omega),\cdots,
m_n(\omega))\right) \\ 
&= (E_D^M)_{\omega}
([m_1(\omega)]_{\omega}^{\circ}
[m_2(\omega)]_{\omega}^{\circ}\cdots
[m_n(\omega)]_{\omega}^{\circ}) 
= 0 
\end{aligned}
$$
for every $(m_1(x),m_2(x),\cdots,m_n(x)) 
\in \Lambda(\{a_k(\omega)\}_{k\in{\mathbb N}}, 
\{b_k(\omega)\}_{k\in{\mathbb N}})$ and 
for every $\omega \in \Omega_0$. 
Since the operation $y \mapsto [y]_{\omega}^{\circ}$ 
is normal and linear, the mapping 
$\Phi_{(A(\omega),B(\omega))}$ 
itself is identically  zero thanks to (c)' 
with the aid of the Kaplansky density theorem. 
Hence, $A(\omega)$ and $B(\omega)$ are free 
with amalgamation over $D(\omega)$ 
for every $\omega \in \Omega_0$ by Lemma 5.2. 
Hence we have proved the assertion. 
\end{proof}  

\medskip\noindent
{\bf 5.2.}\ 
We would like to apply Theorem 5.1 to 
amalgamated free products over Cartan subalgebras. 
In what follows, we suppose that the triple 
$A \supseteq D \subseteq B$ consists of factors 
(with separable preduals) of non-type I and 
a common Cartan subalgebra. The starting point of 
the discussion is Theorem 4.5, (4.6): 
$$
{\mathcal Z}(\widetilde{M}) = 
{\mathcal Z}(\widetilde{A}) \cap {\mathcal Z}(\widetilde{B}) 
\subseteq \widetilde{D}.
$$
Theorem 5.1 implies 

\begin{cor}
Almost every factor $\widetilde{M}(x)$ 
in the central decomposition   
$$
\widetilde{M} = \int_{X_M}^{\oplus} \widetilde{M}(x) 
d\mu(x) 
$$
of the continuous core $\widetilde{M}$ can be 
written as an
amalgamated free product over a common Cartan 
subalgebra.
\end{cor}

Here, we further suppose that $M$ is of type 
III$_{\lambda}$ 
$(0 < \lambda < 1)$, or equivalently that the flow of 
weights 
$(X_M, F_t^M)$ is an (essentially) transitive flow 
with period $-\log\lambda$ (see [CT],[T2]) so that 
we may and do assume that 
$$
X_M = [0,-\log\lambda), \quad 
F_t^M = \text{the translation by $t$}\ (\text{mod}: 
-\log\lambda). 
$$
Since $F_t^M$ is transitive, we have $\widetilde{M}(0) 
\cong 
\widetilde{M}(x)$ for every $x \in X_M = 
[0,-\log\lambda)$, and hence 
we may and do assume that $x \in X_M \mapsto 
\widetilde{M}(x)$ is 
a constant field of the type II$_{\infty}$ factor 
$\widetilde{M}(0)$. 
Hence, we conclude  

\begin{cor}
In the current setting, 
the {\rm (}unique{\rm )} type {\rm II}$_{\infty}$ factor 
appearing in the discrete decomposition {\rm ([C1])} of $M$ 
is written as an amalgamated free product over a common 
Cartan subalgebra. 
\end{cor}

\medskip\noindent
{\bf 5.3.}\ 
We keep the setting and the notations as in \S\S5.2.  
It is known that an ergodic free action of 
a non-amenable discrete group may or may not be amenable 
([Z], see also \S6), or equivalently the associated 
von Neumann factor may or may not be injective {\rm (see [C2])} 
(or hyperfinite), and hence it is somewhat non-trivial 
whether or not the amalgamated free product $M = A *_D B$ 
is non-injective. 

\begin{thm}
In the current setting, 
there is a copy of the free group factor $L({\mathbb F}_2)$ 
in the continuous core $\widetilde{M}$ which is the range 
of a faithful normal conditional expectation. 
In particular, $M$ is not injective. 
\end{thm}

The non-injectivity result follows also from a result in 
our resent work [U2], where we have shown that 
the amalgamated free product $M = A *_D B$ is not a McDuff 
factor (under the assumption that both $A$ and $B$ 
are factors of non-type I). However, the proof below is 
still valid even in the case that both $A$ and $B$ have 
no type I direct summand. (See Remarks 4.8, (2),(3).)

\begin{proof} Let 
$\text{Tr}_{\widetilde{M}}$, 
$\text{Tr}_{\widetilde{A}}$, 
$\text{Tr}_{\widetilde{B}}$, 
$\text{Tr}_{\widetilde{D}}$ 
be the canonical traces on $\widetilde{M}$, 
$\widetilde{A}$, 
$\widetilde{B}$, $\widetilde{D}$, respectively, 
(scaled in the usual way under the dual actions). 
It can be checked that 
$$
\text{Tr}_{\widetilde{M}} 
= \text{Tr}_{\widetilde{D}}\circ\widehat{E_D^M}, \quad 
\text{Tr}_{\widetilde{A}} 
= \text{Tr}_{\widetilde{D}}\circ\widehat{E_D^A}, \quad 
\text{Tr}_{\widetilde{B}} 
= \text{Tr}_{\widetilde{D}}\circ\widehat{E_D^B},  
$$
where $\widehat{E_D^M} : \widetilde{M} \rightarrow 
\widetilde{D}$, 
$\widehat{E_D^A} : \widetilde{A} \rightarrow 
\widetilde{D}$, 
$\widehat{E_D^B} : \widetilde{B} \rightarrow 
\widetilde{D}$ 
are as in \S4. As in \S\S6.1, we consider 
the simultaneous direct integral decompositions of 
the inclusions  
$\widetilde{M} \supseteq \widetilde{A},\ \widetilde{B} 
\supseteq 
\widetilde{D}$ and the conditional expectations 
$\widehat{E^M_D} : \widetilde{M} \rightarrow 
\widetilde{D}$, 
$\widehat{E^A_D} : \widetilde{A} \rightarrow 
\widetilde{D}$, 
$\widehat{E^B_D} : \widetilde{B} \rightarrow 
\widetilde{D}$
subject to the central decomposition 
$$
\widetilde{M} = \int_{X_M}^{\oplus} \widetilde{M}(x)\ 
d\mu(x) 
$$
thanks to Theorem 4.5, (4.6), and let
$$
\text{Tr}_{\widetilde{D}} 
= \int_{X_M}^{\oplus} 
\left(\text{Tr}_{\widetilde{D}}\right)_x\ d\mu(x)
$$
be the corresponding direct integral decomposition. 

\medskip
Note that the continuous core $\widetilde{M}$ is 
of type II$_{\infty}$ or of type II$_1$. (Recall that, 
the continuous core of a von Neumann algebra of type III 
must be of type II$_{\infty}$, while there is no type 
change for the other types.) 
We assume that $\widetilde{M}$ is of type II$_{\infty}$ 
in what follows, since the quite similar (actually  
simpler) argument as below apparently works in the 
type II$_1$ case.

\medskip
Let $X$ be a co-null Borel subset of $X_M$ satisfying 
the same conditions (a), (b), (c) as in the proof of 
Theorem 5.1 with the additional ones: 
\begin{itemize}
\item[(e)] Both $\widetilde{A}(x)$ and $\widetilde{B}(x)$ 
 are von Neumann algebras of type II for every 
 $x \in X$.  (This follows from the assumption 
 that $A$ and $B$ have no type I direct summand.) 
\item[(f)] $\displaystyle{
            \left(\text{Tr}_{\widetilde{D}}\right)_x\circ
            \left(\widehat{E_D^M}\right)_x}$ is a faithful 
            normal semi-finite trace, and thus so are both
            $\displaystyle{
            \left(\text{Tr}_{\widetilde{D}}\right)_x\circ
            \left(\widehat{E_D^A}\right)_x}$, 
            $\displaystyle{
            \left(\text{Tr}_{\widetilde{D}}\right)_x\circ
            \left(\widehat{E_D^B}\right)_x}$. 
\end{itemize}
By choosing a much smaller co-null subset instead of 
$X$ if necessary, the set of those 
$(x,p) \in X\times B({\mathcal H}_0)$ with the separable Hilbert space 
${\mathcal H}_0$ (on which almost every $\widetilde{D}(\omega)$ act),   
satisfying 
\begin{itemize}
\item $p = p^2 = p^* \in \widetilde{D}(x)$; 
\item $\left(\text{Tr}_{\widetilde{D}}\right)_x(p) 
                 = 1$
\end{itemize}
can be assumed to be Borel. Therefore, the measurable 
selection principle enables us to choose a measurable field 
$x \in X \mapsto p_x \in \widetilde{D}(x)$ of projections 
such that $\left(\text{Tr}_{\widetilde{D}}\right)_x(p_x) = 1$ 
for every $x \in X$, and we set 
$$
p := \int_{X_M}^{\oplus} p_x d\mu(x). 
$$
Since $\widetilde{M}$ is of type II$_{\infty}$, 
we have
$$
\begin{aligned}
\int_{X_M}^{\oplus} \widetilde{M}(x)\ d\mu(x) &= 
\widetilde{M} = 
\left(p\widetilde{M}p\right)\otimes B({\mathcal H}_0) \\ 
&= 
\int_{X_M}^{\oplus} 
\left(p_x\widetilde{M}(x)p_x\right)\otimes B({\mathcal H}_0)\ 
d\mu(x) 
\end{aligned}
$$
with the separable Hilbert space ${\mathcal H}_0$. Then we see that, 
for every $x \in X$,  
\begin{itemize}
\item  both $p_x\widetilde{A}(x)p_x$ and
                  $p_x\widetilde{B}(x)p_x$ 
                  are of type II$_1$ (thanks to (e),(f)); 
\item  $p_x\widetilde{A}(x)p_x$, 
                  $p_x\widetilde{B}(x)p_x$ 
                  are free with amalgamation over 
                  $\widetilde{D}(x)p_x$ in the $\widetilde{D}(x)p_x$-probability 
                  space 
                  $\left(p_x\widetilde{M}(x)p_x, 
                  \left(\widehat{E_D^M}\right)_x|_
                  {p_x\widetilde{M}(x)p_x}
                  \right)$ (by (b),(c) and the proof of 
                  Theorem 5.1). 
\end{itemize}
Repeating the argument of the type III$_0$ case 
in the proof of [U1, Lemma 4.2], we choose 
two measurable fields 
$x \in X \mapsto u_x \in p_x\widetilde{A}(x)p_x$, 
$x \in X \mapsto v_x \in p_x\widetilde{B}(x)p_x$
of unitaries satisfying 
$$
\left(\widehat{E_D^A}\right)_x((u_x)^n) = 0, \quad 
\left(\widehat{E_D^B}\right)_x((v_x)^n) = 0
$$
for every $n (\neq 0) \in {\mathbb Z}$, and set 
$$
u := \int_{X_M}^{\oplus} u_x d\mu(x), \quad 
v := \int_{X_M}^{\oplus} v_x d\mu(x).
$$  
Then $u$, $v$ are free Haar unitaries, and hence 
the von Neumann subalgebra $N := \{u, v\}''$ of 
$p\widetilde{M}p$ is 
isomorphic to the free group factor $L({\mathbb F}_2)$. 
Notice that 
$N$ is decomposable relative to ${\mathcal Z}(\widetilde{M}) = 
L^{\infty}(X_M,\mu)$ and that $u_x$, $v_x$ also are 
free Haar unitaries for almost every $x \in X_M$. 
Therefore, we have  
$$
L({\mathbb F}_2) \cong N(x) = \{u_x, v_x\}'' \subseteq 
p_x\widetilde{M}(x)p_x \quad 
\text{for almost every $x \in X_M$}, 
$$
and thus 
$$
\begin{aligned}
L({\mathbb F}_2)\otimes B({\mathcal H}_0)\otimes{\mathbf C}1 
&\subseteq 
L({\mathbb F}_2)\otimes B({\mathcal H}_0)\otimes{\mathcal Z}(\widetilde{M}) \\ 
&= 
\int_{X_M}^{\oplus} L({\mathbb F}_2)\otimes B({\mathcal H}_0)\ d\mu(x) \\ 
&\subseteq 
\int_{X_M}^{\oplus} \left(p_x\widetilde{M}(x)p_x\right)
\otimes B({\mathcal H})\ d\mu(x) \\
&= 
\int_{X_M}^{\oplus} \widetilde{M}(x)\ d\mu(x). 
\end{aligned} 
$$
Since $p\widetilde{M}p$ is of type II$_1$, the copy of 
$L({\mathbb F}_2)\otimes B({\mathcal H}_0)$ in 
$\widetilde{M}$ (or in $\widetilde{M}(x)$ for almost every 
$x \in X_M$) is clearly the range of a faithful normal conditional 
expectation. Hence we are done. 
\end{proof}

\medskip
Since the copy of $L({\mathbb F}_2)$ constructed in the 
proof is well-behaved with the central decomposition 
of $\widetilde{M}$, we have 

\begin{cor}
Keep the same setting as in Theorem 5.5.  
If the amalgamated free product $M = A *_D B$ is of 
type {\rm III}$_{\lambda}$, then the type {\rm II}$_{\infty}
$ factor appearing in the discrete decomposition of 
$M$ contains a copy of the free group factor 
$L({\mathbb F}_2)$ which is the range of a faithful normal 
conditional expectation. Therefore, so does the $M$ itself. 
\end{cor} 
\begin{proof} The first part of the assertion is clear 
from the proof of Theorem 5.5 together with Corollary 5.4. 
The latter follows from the first half and the discrete 
decomposition theorem ([C1]) for type III$_{\lambda}$ factors. 
\end{proof}  

\medskip
For a while, we assume that both $A$ and $B$ are 
general factors (not necessary of non-type I) and that 
$D$ is a common Cartan subalgebra. If either $A$ or $B$ 
is of type I$_n$ with possibly $n = \infty$, then both 
must coincide, i.e., $A = B =  M_n({\mathbf C})$ or 
$B({\mathcal H})$. In this case, we can see that 
$M = A  *_D B$ ($= A *_D A$) is isomorphic to 
$L({\mathbb F}_{n-1})\otimes A$. Therefore, we obtain 

\begin{cor}
The amalgamated free product $M = A *_D B$ of factors 
{\rm (}with separable preduals{\rm )} over a common Cartan 
subalgebra is injective if and only if  either $A$ or $B$ is 
of type {\rm I}$_2$ {\rm (}and hence $A = B = M_2({\mathbf C})$ 
and $D$ is the diagonals{\rm)}. 
\end{cor}

This corollary can be thought of as an analogue of 
the following classical group theoretical fact: 
A free product group $G * H$ is amenable if and only if 
$G = H = {\mathbb Z}_2$. 

\medskip\noindent
{\bf 5.4.}\ We keep the same setting and the notations 
as in \S\S5.2 even in this subsection. We would like here 
to show that the amalgamated free  product $M = A *_D B$ 
is not related to any free group factor as a simple 
application of the striking result [V2] of D. Voiculescu 
with the aid of Theorem 5.1 (or Corollary 5.3, Corollary 5.4). 
A similar application of Voiculescu's result was also given 
by D. Shlyakhtenko [S1] in a different context. 

\medskip
When the amalgamated free product $M = A *_D B$ is 
of type II$_1$, we can apply directly Voiculescu's theorem 
to the case since the  normalizer ${\mathcal N}_M(D)$ generates 
the whole $M$, and hence $M$ is not isomorphic to 
any (interpolated) free group factor $L({\mathbb F}_r)$ 
with $0 < r \leq  \infty$. Thus it suffices to consider 
only the infinite cases, and we start with 
the following lemma: 

\begin{lemma}
Let $N \supseteq C$ be a factor of type {\rm II}$_{\infty}$ 
and an abelian von Neumann subalgebra. 
Let $\text{\rm Tr}$ be a faithful normal semi-finite 
trace on $N$ such that $\text{\rm Tr}|_C$ is semi-finite. 
Suppose that the normalizer ${\mathcal N}_N(C)$ generates 
the whole $N$. Then, for each finite {\rm (}in $N${\rm )} 
non-zero projection $p \in C$ {\rm (}such a projection 
indeed exists since $\text{\rm Tr}|_C$ is semi-finite{\rm )}, 
the normalizer ${\mathcal N}_{pNp}(Cp)$ generates the whole 
$pNp$. 
\end{lemma}
\begin{proof} By assumption, we see that the linear 
span of ${\mathcal N}_N(C)$ forms a $\sigma$-weakly dense 
$*$-subalgebra in $N$. Hence, the linear span of elements 
of the form $pup$ with $u \in {\mathcal N}_N(C)$ also is 
$\sigma$-weakly dense in $pNp$. 
Since $p$ is finite in $N$, $pNp$ is a factor 
of type $II_1$. Let us denote $q$ and $r$ the support and 
the range projections of $pup$, respectively. 
Then the projections $q, r$ are in $pNp$ since $q, r \leq p$. 
We here need the following fact: 

\medskip
\begin{fact}
If $p \neq q$ {\rm (}or equivalently $p \neq r$ 
thanks to the fact that $pNp$ is finite{\rm )}, 
then there is an element $w \in {\mathcal G}{\mathcal N}_{pNp}(Cp)$ 
with $w^* w = p-q$, $ww^* = p-r$. 
Here, ${\mathcal G}{\mathcal N}_{pNp}(Cp)$ denotes the normalizing 
groupoid, i.e., the set of  those partial isometries 
$v \in N$  such that $v^* v, vv^* \in C$ and 
$v C v^* = C vv^*$, 
$v^* C v = Cv^* v$. 
\end{fact} 
\begin{proof} (Proof of the Fact.)\   
Since $N$ is a factor, we have $(p-r)N(p-q) \neq \{0\}$, 
and hence there is a unitary $v \in {\mathcal N}_N(C)$ 
such that $(p-r)v(p-q)$ is not equal zero. Thus, there is 
a non-zero element $v_0 \in {\mathcal G}{\mathcal N}_{pNp}(Cp)$ 
such that $v_0^* v_0 \leq  p-q$ and $v_0 v_0^* \leq p-r$. 
We can do the standard exhaustion argument thanks 
to the fact that $pNp$ is finite. Hence we get 
a desired partial isometry. 
\end{proof}

Let $w := pup$ with $u \in {\mathcal N}_N(C)$, and 
the above fact says that we can choose a unitary 
$\widetilde{w} \in {\mathcal N}_{pNp}(Cp)$ in such a way that 
$w = r\widetilde{w}q$. Moreover, $q, r$ can be written as 
finite linear combinations of unitries in $Cp$, and 
hence $r\widetilde{w}q$ is a finite linear combinations 
of elements in ${\mathcal N}_{pNp}(Cp)$. 
Therefore, any element in $pNp$ can be approximated 
$\sigma$-weakly by finite linear combinations of 
elements in ${\mathcal N}_{pNp}(Cp)$. Hence we complete 
the proof of the lemma. 
\end{proof}

\begin{prop}
Keep the same setting as in Theorem 5.5.  

\noindent
{\rm (1)} 
If the amalgamated free product $M = A *_D B$  
is of type {\rm II}$_{\infty}$, then $M$ is not 
isomorphic to any $L({\mathbb F}_r)\otimes B({\mathcal H})$ 
with $0 < r \leq \infty$. 

\noindent
{\rm (2)} 
If the amalgamated free product $M = A *_D B$ 
is of type {\rm III}, then almost every 
type {\rm II}$_{\infty}$ factor
appearing in  the central decomposition
of the continuous  core
$\widetilde{M}$ is not isomorphic to any 
$L({\mathbb F}_r)\otimes B({\mathcal H})$ with 
$0 < r \leq \infty$. 

\noindent
{\rm (3)} 
If the amalgamated free product $M = A *_D B$ 
is of type {\rm III}$_{\lambda}$ 
{\rm (}$0 < \lambda < 1${\rm)}, 
then the type {\rm II}$_{\infty}$ factor 
appearing in the discrete decomposition 
is not isomorphic to 
any $L({\mathbb F}_r)\otimes B({\mathcal H})$ with 
$0 < r \leq \infty$. 
\end{prop}
\begin{proof} All the assertions follow from 
[V2, 5,3 Theorem, 7.4 Corollary] with the aid 
of Lemma 5.8. When showing the assertions 
(2), (3), we further need Corollary 5.3, 
Corollary 5.4, respectively. 
\end{proof}

\section{Example I.\ \ Boundary Actions of Free Groups}

\medskip\noindent
{\bf 6.1.}\ 
Let ${\mathfrak X}$ be a finite set with $|{\mathfrak X}| \geq 2$, 
and we set ${\mathfrak X}^{-1} :=  \{x^{-1} ; x \in {\mathfrak X}\}$. 
We consider the free group $\Gamma := 
{\mathbb F}({\mathfrak X})$ over the generators ${\mathfrak X}$ 
and its boundary $\partial\Gamma$. In this case, 
the boundary $\partial\Gamma$ is defined as 
the one-sided shift space of the alphabets 
${\mathfrak X}\cup{\mathfrak X}^{-1}$ determined by 
the forbidden blocks $(x^{-1} x)$, $(x x^{-1})$, 
equipped with the usual product topology. 
It is plain to see that $\partial\Gamma$ is identified 
with the set of semi-infinite reduced words in 
${\mathfrak X}\cup{\mathfrak X}^{-1}$. We will freely use 
these two different descriptions in what follows. 
The group $\Gamma$ acts topologically on the boundary 
$\partial\Gamma$ by the left multiplication, i.e., 
for $\gamma \in \Gamma$ and for 
$\omega = \omega_1\omega_2\cdots \in \partial\Gamma$,  
$$
\gamma\cdot\omega := 
\text{the reduced form of the word 
$\gamma\omega_1\omega_2\cdots$}. 
$$

\medskip\noindent
{\bf 6.2.}\ 
We decompose the set ${\mathfrak X}$ into two disjoint 
non-empty subsets ${\mathfrak X}_1$, ${\mathfrak X}_2$ 
with ${\mathfrak X} = {\mathfrak X}_1 \sqcup {\mathfrak X}_2$. 
Then we have $\Gamma = \Gamma_1 * \Gamma_2$ 
with $\Gamma_1 = {\mathbb F}({\mathfrak X}_1)$, $\Gamma_2 = 
{\mathbb F}({\mathfrak X}_2)$.  
The boundaries 
$\partial\Gamma_1$, 
$\partial\Gamma_2$ 
can be (topologically) embedded into $\partial\Gamma$ 
as follows: 
$$
\begin{aligned}
\partial\Gamma_1 
= &\text{the shift subspace of $\partial\Gamma$} \\
  &\text{with extra forbidden blocks $(x)$, $x \in {\mathfrak X}_2$},\\
\partial\Gamma_2 
= &\text{the subset of $\partial\Gamma$} \\
  &\text{with extra forbidden blocks $(x)$, $x \in {\mathfrak X}_1$}. 
\end{aligned}
$$
Therefore, the subspaces $\partial\Gamma_1$, 
$\partial\Gamma_2$ are closed and invariant under 
the actions of $\Gamma_1$, $\Gamma_2$, respectively. 
We consider the following disjoint decompositions: 
$$
\partial\Gamma = 
(\partial\Gamma_1)^c \sqcup \partial\Gamma_1, \quad 
\partial\Gamma = 
(\partial\Gamma_2)^c \sqcup \partial\Gamma_2, 
$$
and note that $(\partial\Gamma_1)^c$, 
$(\partial\Gamma_2)^c$ are open and invariant 
under the actions of $\Gamma_1$, $\Gamma_2$, 
respectively. We define 
$$
\begin{aligned}
(\partial\Gamma_1)^{\perp} &:= 
\{ \omega = (\omega_n)_{n = 1}^{\infty} \in 
(\partial\Gamma)^c\ ;\ 
\omega_1 \in {\mathfrak X}_2 \cup ({\mathfrak X}_2)^{-1} \}, \\
(\partial\Gamma_2)^{\perp} &:= 
\{ \omega = (\omega_n)_{n = 1}^{\infty} \in 
(\partial\Gamma)^c\ ;\ 
\omega_1 \in {\mathfrak X}_1 \cup ({\mathfrak X}_1)^{-1} \}. 
\end{aligned}
$$
Let us define the map 
$\Phi_1 : \Gamma_1\times(\partial\Gamma_1)^{\perp} 
\rightarrow (\partial\Gamma_1)^c$ 
as the restriction of the action map 
$\Gamma \times \partial\Gamma \rightarrow 
\partial\Gamma$, i.e., 
\begin{equation}
\Phi_1\ :\ 
\begin{cases}
(e, \omega) \mapsto \omega, \\
(\gamma_1\cdots\gamma_2, (\omega_n)_{n = 1}^{\infty}) 
\mapsto 
(\gamma_1,\dots,\gamma_n,\omega_1,\omega_2,\dots) 
\end{cases} \tag{6.2.1}
\end{equation}
for each reduced word $\gamma_1\cdots\gamma_n$. 
Similarly, the map 
$\Phi_2 : \Gamma_2\times(\partial\Gamma_2)^{\perp} 
\rightarrow (\partial\Gamma_2)^c$ is defined as 
the restriction of the action map
$\Gamma \times \partial\Gamma \rightarrow \partial\Gamma$. 
We here note that the topology on $\partial\Gamma$ is 
generated by the family of clopen sets of the form: 
$$
\Omega(\gamma) = 
\{ \omega = (\omega_n)_{n = 1}^{\infty} \in 
\partial\Gamma\ ;\ 
\omega_1 = \gamma_1,\dots,\omega_n = \gamma_n \}
$$
with a reduced word $\gamma = \gamma_1\cdots\gamma_n$. 

\begin{lemma}
We have
\begin{gather}
(\partial\Gamma_1)^c = 
\bigsqcup_{\gamma \in \Gamma_1} 
\bigsqcup_{\omega \in {\mathfrak X}_2 \cup 
({\mathfrak X}_2)^{-1}} 
\Omega(\gamma \omega), \tag{6.2.2}
\\
(\partial\Gamma_1)^c = 
\bigsqcup_{\gamma \in \Gamma_2} 
\bigsqcup_{\omega \in {\mathfrak X}_1 \cup 
({\mathfrak X}_1)^{-1}} 
\Omega(\gamma\omega). \tag{6.2.3}
\end{gather}
\end{lemma}

\medskip
It is plain to check that  
\begin{enumerate}
\item the maps $\Phi_1$, $\Phi_2$ are bijections; 
\item $
\Phi_k(\{\gamma\}\times\Omega(\omega_1\cdots \omega_n)) 
= 
\Omega(\gamma \omega_1 \cdots \omega_n)$ ($k = 1,2$) 
\newline 
with a reduced word $\omega_1\cdots \omega_n$. 
\end{enumerate}
Therefore, thanks to Lemma 6.1, we see that $\Phi_1$, 
$\Phi_2$ send 
all the basic clopen sets to all those, when the 
product topologies of 
the discrete one and the induced one from 
$\partial\Gamma$ are considered 
on both 
$\Gamma_1\times(\partial\Gamma_1)^{\perp}$, 
$\Gamma_2\times(\partial\Gamma_2)^{\perp}$. 
It is also plain to check that 
$$
\Phi_k(\gamma_1\gamma_2, \omega) = 
\gamma_1\cdot\Phi_k(\gamma_2,\omega)  
\quad (k = 1,2), 
$$
and hence we conclude 

\begin{prop}
The maps 
$\Phi_1 : \Gamma_1\times(\partial\Gamma_1)^{\perp} 
\rightarrow (\partial\Gamma_1)^c$, 
$\Phi_2 : \Gamma_2\times(\partial\Gamma_2)^{\perp} 
\rightarrow (\partial\Gamma_2)^c$
are homeomorphisms, and via these homeomorphisms, the 
actions 
of $\Gamma_1$, $\Gamma_2$ on $(\partial\Gamma_1)^c$, 
$(\partial\Gamma_2)^c$ are conjugate to those of 
$\Gamma_1$, $\Gamma_2$ 
on $\Gamma_1\times (\partial\Gamma_1)^{\perp}$,
$\Gamma_1\times(\partial\Gamma_2)^{\perp}$ which 
are defined as the product action of the translation 
and the trivial one.  
\end{prop}  

\medskip\noindent
{\bf 6.3.}\ 
Let $n = |{\mathfrak X}|$ and $n_1 = |{\mathfrak X}_1|$, $n_2 = 
|{\mathfrak X}_2|$. 
We discuss here the probability measure $\mu$ 
on the boundary $\partial\Gamma$ defined in such a way 
that  
$$
\mu(\Omega(\gamma)) := 
\frac{1}{2n}\left(\frac{1}{2n - 1}\right)^{\ell(\gamma) 
- 1}
$$
with the word length function $\ell(\ \cdot\ )$. 
It is known that the measure $\mu$ is quasi-invariant 
under the action of $\Gamma$. (See [KS].) 
The non-singular action of $\Gamma$ on the probability 
space $(\partial\Gamma, \mu)$ can be checked to be free 
and ergodic (see [RR],[KS],[PS]). 
Moreover, J. Ramagge \& G. Robertson [RR] showed that 
the action of $\Gamma$ is of type III$_{\frac{1}{2n-1}}$ 
so that the crossed-product 
$M = L^{\infty}(\partial\Gamma,\mu)\rtimes\Gamma$ 
is a factor of type III$_{\frac{1}{2n-1}}$. 
Moreover, S. Adams' result [A] (see also [Ver, Example 2 in p\. 89]) 
implies that the factor is injective. (It should be remarked that 
the Cuntz-Krieger algebra interpretation for boundary actions 
provided by J. Spielberg [Sp] together with M. Enomoto, M. Fujii 
\& Y. Watatani [EFW] also shows that the crossed-product is 
the injective factor of type III$_{\frac{1}{2n-1}}$.) Set 
$A := L^{\infty}(\partial\Gamma,\mu)\rtimes\Gamma_1$, 
$B := L^{\infty}(\partial\Gamma,\mu)\rtimes\Gamma_2$, 
and it is plain to see that, the crossed-product is 
written as an amalgamated free product over a common 
Cartan subalgebra, that is, we have $M \cong A *_D B$ 
with $D := L^{\infty}(\partial\Gamma,\mu)$. 
Therefore, the boundary action of the free group $\Gamma$ 
provides an example of an injective factor arising as 
an amalgamated free product over a common Cartan subalgebra. 

\medskip\noindent
{\bf 6.4.}\ Since  
$$
|\{ \gamma \in \Gamma_k\ ;\ \ell(\gamma) = m \}| = 
2n_k (2n_k - 1)^{m - 1} \quad (k = 1,2), 
$$
we have  
$$
\begin{aligned}
&\mu((\partial\Gamma_1)^c) = 
\mu\left(
\bigsqcup_{\gamma \in \Gamma_1} 
\bigsqcup_{\omega \in {\mathfrak X}_2 \cup ({\mathfrak 
X}_2)^{-1}} 
\Omega(\gamma \omega)
\right) 
\quad 
\text{(Lemma 6.1)} \\
&= 
\sum_{\omega \in {\mathfrak X}_2 \cup ({\mathfrak X}_2)^{-1}}
\mu(\Omega(\omega)) + 
\sum_{m = 1}^{\infty}
\sum_{\scriptstyle \gamma \in \Gamma_1 \atop \scriptstyle \ell(\gamma) = m} 
\sum_{\omega \in {\mathfrak X}_2 \cup ({\mathfrak X}_2)^{-1}}
\mu(\Omega(\gamma\omega)) \\
&= 
2n_2 \frac{1}{2n} + 
\sum_{m = 1}^{\infty} 
\left(2n_1(2n_1 - 1)^{m - 1}\right) 2n_2 
\left(
\frac{1}{2n}
\left(\frac{1}{2n - 1}\right)^m
\right) \\
&= 
\frac{n_2}{n} + 
\frac{2 n_1 n_2}{n(2n-1)} \cdot 
\sum_{m = 1}^{\infty} 
\left(\frac{2n_1 - 1}{2n - 1}\right)^{m - 1} \\
&= 1 
\quad \text{(by $n = n_1 + n_2$)}.
\end{aligned}
$$
Similarly, we have $\mu((\partial\Gamma_2)^c) = 1$. 
Hence, we obtain  
\begin{equation}
\mu(\partial\Gamma_1) = \mu(\partial\Gamma_2) = 0. \tag{6.4.1}  
\end{equation}
Furthermore, we have 
$$
\begin{aligned}
\mu(
\Phi_k(\{\gamma\}\times&\Omega(\omega_1\cdots \omega_m))
) = 
\frac{1}{2n}
\left(\frac{1}{2n-1}\right)^{\ell(\gamma) + m - 1} \\
&= 
\left(\frac{1}{2n-1}\right)^{\ell(\gamma)} \times 
\frac{1}{2n}\left(\frac{1}{2n-1}\right)^{m-1}, 
\end{aligned}
$$
and hence 
\begin{equation}
(\mu|_{(\partial\Gamma_k)^c})\circ\Phi_k = 
\delta_k\otimes(\mu|_{(\partial\Gamma_k)^{\perp}}) 
\quad (k = 1,2) \tag{6.4.2} 
\end{equation}
with the measure $\displaystyle{\delta_k(\{\gamma\}) = 
\left(\frac{1}{2n - 1}\right)^{\ell(\gamma)}}$ 
equivalent to 
the counting measure. From the discussions above, we 
conclude 

\begin{prop}
We have
\begin{align}
\begin{split}
&L^{\infty}(\partial\Gamma, \mu)\rtimes\Gamma_1 = 
L^{\infty}((\partial\Gamma_1)^c, \mu)\rtimes\Gamma_1 
\\
&\cong 
\left(\ell^{\infty}(\Gamma_1)\rtimes\Gamma_1\right)
\otimes L^{\infty}((\partial\Gamma_1)^{\perp}, 
\mu|_{(\partial\Gamma_1)^{\perp}}), 
\end{split} \tag{6.4.3}\\
\begin{split}
&L^{\infty}(\partial\Gamma, \mu)\rtimes\Gamma_2 = 
L^{\infty}((\partial\Gamma_2)^c, \mu)\rtimes\Gamma_2 \\
&\cong 
\left(\ell^{\infty}(\Gamma_2)\rtimes\Gamma_2\right)
\otimes L^{\infty}((\partial\Gamma_2)^{\perp}, 
\mu|_{(\partial\Gamma_2)^{\perp}}). 
\end{split}
\tag{6.4.4} 
\end{align}
The isomorphisms are induced from the maps $\Phi_1$, 
$\Phi_2$, respectively. 
In particular, the crossed-products both are of 
homogeneous 
type {\rm I}$_{\infty}$. 
\end{prop}

Therefore, we have seen that the free components 
of our injective amalgamated free product $M = A *_D B$ 
both are of homogeneous type I$_{\infty}$. 

\bigskip\noindent
{\bf 6.5.}\ 
At the end of this section, we give a criterion on 
injectivity of amalgamated free products over Cartan 
subalgebras. Let $M = A *_D B$ be an amalgamated free 
product over a common Cartan subalgebra. Here, we do not 
assume that $A$ and $B$ are factors nor that they have 
no type I direct summand. We choose central projections 
$p_A$, $p_B$ of $A$, $B$ in such a way that both $Ap_A$ 
and $Bp_B$ have no type I direct summand. Since $D$ is 
a Cartan subalgebra in both $A$ and $B$, the projections 
$p_A$, $p_B$ are in $D$ so that $p := p_A p_B = p_B p_A$ 
is also a projection in $D$. Suppose here that $p$ is 
non-zero. Then the reduced von Neumann algebra $pMp$ contains 
both $pAp$ and $pBp$, and their freeness can be easily checked 
with respect to the conditional expectation 
$$
(E_D^M)_p := E_D^M|_{pMp} : pMp \rightarrow Dp.  
$$
It can be easily checked that the continuous cores satisfy  
\begin{equation}
\widetilde{pMp} = p\widetilde{M}p, \quad 
\widetilde{pAp} = p\widetilde{A}p, \quad 
\widetilde{pBp} = p\widetilde{B}p, \quad 
\widetilde{Dp} =  \widetilde{D}p. 
\tag{6.5.1}
\end{equation}  
Hence we get the inclusion relations 
\begin{equation}
\widetilde{pMp} \supseteq 
\widetilde{pAp} \supseteq 
\widetilde{Dp}, \quad 
\widetilde{pMp} \supseteq 
\widetilde{pBp} \supseteq 
\widetilde{Dp}. 
\tag{6.5.2} 
\end{equation}
We can easily see that the conditional expectation 
$$
\widetilde{(E_D^M)_p} : 
\widetilde{pMp} \rightarrow \widetilde{Dp}
$$ 
coincides with 
$$
\left(\widetilde{E_D^M}\right)_p := 
\widetilde{E_D^M}|_{p\widetilde{M}p} : 
p\widetilde{M}p \rightarrow \widetilde{D}p.
$$
Hence we can show that, the von Neumann subalgebra 
$$
N := 
\widetilde{pAp} \vee \widetilde{pBp}\ 
(\subseteq \widetilde{pMp})
$$
is identified with the amalgamated free product over 
a common Cartan subalgebra 
$$
\left(
\widetilde{pAp}
\right) 
*_{\widetilde{Dp}} 
\left(
\widetilde{pBp}
\right). 
$$
Notice here that 
\begin{itemize} 
\item both $pAp$ and $pBp$ have no type I direct summand; 
\item there is a faithful normal conditional expectation 
from $\widetilde{pMp}$ onto $N$ since $N$ is invariant 
under the modular action 
$\sigma_t^{\psi\circ\widetilde{(E_D^M)_p}}$ 
($t \in {\mathbf R}$) with a faithful normal state $\psi$ 
on $\widetilde{Dp}$ thanks to [T1] and Theorem 1.1. 
\end{itemize}
Thus we apply the same argument as in the proof of
Theorem 5.5  (see after the statement of that theorem) to 
the amalgamated free product $N$, and as a consequence 
we get a copy of the free group factor in $\widetilde{pMp}$ 
as the range of a faithful normal conditional expectation 
from $N$ (and hence from $\widetilde{pMp}$). Therefore, 
$pMp$ is not injective, and neither is $M$. Therefore, 
we conclude

\begin{prop}
In the current setting, 
if the amalgamated free product $M = A *_D B$ is injective, 
then the non-type I direct summands in $A$ and $B$ need not 
meet {\rm (}in $D${\rm )}, i.e., their support central 
projections are disjoint {\rm (}in $D${\rm )}. 
\end{prop}

\noindent 
\begin{remark} {\rm One can construct two von Neumann algebras 
$A$, $B$ with a common Cartan subalgebra $D$ in such a way 
that (i) $A$ has the non-type I direct summand and $B \cong 
L^{\infty}(\Omega)\otimes B({\mathcal H})$ (possibly with any 
dimension $\text{dim}{\mathcal H} \geq 2$); (ii) the amalgamated 
free product $A *_D B$ is injective. 
(Compare with Proposition 6.3, 
6.4.)}   
\end{remark}

\section{Example II.\ \  Number of Free Components}

Let $A \supseteq D \subseteq B$ be $\sigma$-finite von 
Neumann 
algebras with faithful normal conditional expectations 
$E_D^A : A \rightarrow D$, $E_D^B : B \rightarrow D$, 
which  
are assumed to be of the form: 
$$
A = A_0\otimes B(\ell^2({\mathbb N})), \quad 
B = B_0\otimes B(\ell^2({\mathbb N})), \quad 
D = D_0\otimes \ell^{\infty}({\mathbb N}), 
$$
$$
E_D^A(a\otimes e_{ij}) = \delta_{ij} E_{D_0}^{A_0}(a), 
\quad 
E_D^A(b\otimes e_{ij}) = \delta_{ij} E_{D_0}^{B_0}(b)
$$
with 
$$
B(\ell^2({\mathbb N})) = \{e_{ij}\}'' \supseteq 
\ell^{\infty}({\mathbb N}) 
= \{ e_{ii} \}'',
$$ 
where the $e_{ij}$'s are the natural matrix units and 
will be denoted by $e_{ij}^A$ or $e_{ij}^B$ instead of 
$e_{ij}$ when regarded as elements in $A$ or $B$, 
to avoid any confusion.  We further suppose that $B_0$ 
(and hence $B$ itself) is injective or hyperfinite, has 
no type I direct summand, and that $D_0$ is a Cartan 
subalgebra. Thanks to A. Connes, J. Feldman \& B. Weiss 
[CFW], we may and do assume that there is a unitary 
$u \in B_0$ such that 
$$
B_0 = <D_0, u >'', \quad E_{D_0}^{B_0}(u^n) = 0\ 
\text{as long as $n \neq 0$},
$$
and the automorphism $\text{Ad}u \in \text{Aut}(D_0)$ 
is denoted by $\alpha$. 

In this setting, we will investigate the reduced 
von Neumann algebra $pMp$ of the amalgamated free 
product: 
$$
(M, E_D^M) = (A, E_D^A) *_D (B, E_D^B) 
$$ 
by a minimal projection $p := 1\otimes e_{11}$ 
in the common subalgebra 
${\mathbf C}1\otimes \ell^{\infty}({\mathbb N})$.  

\medskip
We introduce the following notation rule: 
$$
[a]_{ij}^A := a\otimes e_{ij}^A \quad \text{in}\ A, 
\quad \text{and} \quad  
[b]_{ij}^B := b\otimes e_{ij}^B \quad \text{in}\ B, 
$$
and, in what follows, will freely use 
the identification: 
$$
e_{ij}^A = 1\otimes e_{ij}^A = [1_{A_0}]_{ij}^A, \quad 
e_{ij}^B = 1\otimes e_{ij}^B = [1_{B_0}]_{ij}^B. 
$$

\begin{lemma}
We have  
\begin{gather}
[a_1]_{ij}^A [a_2]_{k\ell}^A = 
\delta_{jk} [a_1 a_2]_{i\ell}^A, \quad     
[b_1]_{ij}^B [b_2]_{k\ell}^B = 
\delta_{jk} [b_1 b_2]_{i\ell}^B; 
\tag{7.1}\\
\begin{split}
e_{ij}^A [a]_{k\ell}^A = 
\delta_{jk} [a]_{i\ell}^A 
= [a]_{ij}^A e_{k\ell}^A, \\
e_{ij}^B [b]_{k\ell}^B = 
\delta_{jk} [b]_{i\ell}^B 
= [b]_{ij}^B e_{k\ell}^B; 
\end{split}\tag{7.2}\\
[a]_{ij}^A{}^* = [a^*]_{ji}^A, \quad 
[b]_{ij}^B{}^* = [b^*]_{ji}^B. 
\tag{7.3}
\end{gather}
\end{lemma}

\begin{lemma}
The von Neumann algebra $B$ is 
generated by $D_0\otimes{\mathbf C}1$ and partial 
isometries $[u^n]_{i1}^B$, $n \in {\mathbb Z}$, 
$i  = 2,3\dots$  
\end{lemma}

\medskip
Set $u(n,i) := e_{1i}^A\cdot[u^n]_{i1}^B$, 
$n \in {\mathbb Z}$, $i = 2,3,\dots$, a unitary in $pMp$.  

\begin{lemma}
We have, for each $n \in {\mathbb Z}$, 
$i = 2,3,\dots$,  
$$
E_D^M(u(n,i)^k) = 0
$$ 
whenever $k \neq 0\ (\in {\mathbb Z})$. 
\end{lemma}
\begin{proof} One may and do assume $k > 0$ since 
$E_D^M(u(n,i)^{-k}) = E_D^M(u(n,i)^k)^*$. Notice that  
$$
\begin{cases}  
E_D^M(e_{1i}^A) = E_D^A(e_{1i}^A) = 0 
&\text{as long as}\ i \neq 1, \\ 
E_D^M([u^n]_{i1}^B) = E_D^B(u^n\otimes e_{i1}^B) 
= \delta_{i1}\cdot E_{D_0}^{B_0}(u^n) = 0 
&\text{as long as}\ i \neq 1.   
\end{cases}
$$
Therefore, by the freeness, we have, for $n \in 
{\mathbb Z}$, 
$i \neq 1$,  
$$
E_D^M(u(n,i)^k) = 
E_D^M(e_{1i}^A\cdot[u^n]_{i1}^B\cdots 
e_{1i}^A\cdot[u^n]_{i1}^B) = 0. 
$$
Hence we are done. 
\end{proof}

We define the faithful normal conditional expectation 
$$
(E_D^M)_p := E_D^M|_{pMp} : 
pMp \rightarrow Dp = D_0\otimes{\mathbf C}p
$$ 
(which is well-defined since $p$ is in 
the smaller algebra $D$). 

\begin{lemma}
The family 
$$
\{ A_0\otimes{\mathbf C}p \} \cup 
\{ u(n,i) : n \in {\mathbb Z}, i = 2,3,\dots \}
$$
is free with amalgamation over $D_0\otimes{\mathbf C}p$ 
with respect to 
$(E_D^M)_p$. 
\end{lemma}
\begin{proof} Since all the $u(n,i)$'s normalize 
the subalgebra $D_0\otimes{\mathbf C}p$ and since 
$(D_0\otimes{\mathbf C}p)(A_0^{\circ}\otimes{\mathbf C}p)$ 
and 
$(A_0^{\circ}\otimes{\mathbf C}1)(D_0\otimes{\mathbf C}1)$ 
are contained in $(A_0^{\circ}\otimes{\mathbf C}p)$, 
it suffices to show that 
\begin{equation}
\begin{split}
E_D^M([a_1]_{11}^A &u(n_1,i_1)^{k_1} [a_2]_{11}^A u(n_2, i_2)^{k_2} \\
&\cdots[a_m]_{11}^A u(n_m, i_m)^{k_m} [a_{m+1}]_{11}^A) = 0 
\end{split} \tag{7.4}
\end{equation}
whenever all $k_j$'s are not equal to $0$, 
the beginning and the ending letters $a_1$,$a_{m+1}$ are 
the identity $1$ or in $A_0^{\circ} 
:= \text{Ker} E_{D_0}^{A_0}$, and the other $a_j$'s 
are 
\begin{equation}
a_j\ \text{is}\  
\begin{cases}
\text{in $A^{\circ}$ or the identity $1$} & 
\text{if}\ (n_{j-1}, i_{j-1}) \neq (n_j, i_j), \\
\text{in}\ A_0^{\circ} & 
\text{if}\ (n_{j-1}, i_{j-1}) = (n_j, i_j).
\end{cases}  \tag{7.5}
\end{equation}

We have, for $a \in A$,  
$$
\begin{aligned}
u&(n_1,i_1)^{-k_1}\cdot[a]_{11}^A u(n_2,i_2)^{k_2} 
\\
&= 
[u^{-n_1}]_{1 i_1}^B e_{i_1 1}^A
\cdots
[u^{-n_1}]_{1 i_1}^B [a]_{i_1 i_2}^A 
[u^{n_2}]_{i_2 1}^B
\cdots
e_{1 i_2}^A\cdot[u^{n_2}]_{i_2 1}^B, \\
u&(n_1,i_1)^{k_1} [a]_{11} u(n_2,i_2)^{-k_2} 
\\
&= 
e_{1 i_1}^A [u^{n_1}]_{i_1 1}^B\cdots 
e_{1 i_1}^A [u^{n_1}]_{i_1 1}^B [a]_{11}^A 
[u^{-n_2}]_{1 i_2}^B e_{i_2 1}^A
\cdots[u^{-n_2}]_{1 i_2}^B e_{i_2 1}^A, 
\end{aligned}
$$
Thus, if $k_1, k_2 \neq 0$ and if $a$ is such as in 
(7.5), we see that 
$$
\begin{aligned}
u(n_1,i_1)^{-k_1} [a]_{11}^A u(n_2,i_2)^{k_2} 
&\in 
[u^{-n_1}]_{1 i_1}^B 
\underbrace{A^{\circ}\cdots A^{\circ}}_{alternating} 
[u^{n_2}]_{i_2 1}^B \\
&\subseteq 
B^{\circ} A^{\circ}\cdots A^{\circ} B^{\circ}, \\ 
u(n_1,i_1)^{k_1} [a]_{11} u(n_2,i_2)^{-k_2} 
&\in  
u(n_1,i_1) 
\underbrace{A^{\circ}\cdots A^{\circ}}_{alternating} 
u(n_2,i_2)^* \\ 
&\subseteq 
A^{\circ}\cdot B^{\circ}\cdot A^{\circ} \cdots
A^{\circ}\cdot B^{\circ}\cdot A^{\circ}. 
\end{aligned}
$$
Notice that 
\begin{gather}
[u^{n_1}]_{1 i_1}^B [u^{-n_2}]_{1 i_2}^B 
= [u^{n_1 - n_2}]_{i_1 i_2}^B \in B^{\circ}, \notag\\
u(n_1, i_1)^* u(n_2, i_2) 
= \begin{cases} 
[u^{-n_1}]_{1 i_1}^B e_{i_1 i_2}^A 
[u^{n_2}]_{i_2 1}^B \in 
B^{\circ} A^{\circ} B^{\circ} &(i_1 \neq i_2), \\
[u^{n_2 - n_1}]_{1 1}^B \in B^{\circ} &(i_1 = i_2)  
\end{cases}\notag
\end{gather}
as long as $(n_1,i_1) \neq (n_2,i_2)$, and 
one can easily check the desired equality (7.4) 
based on the above facts. 
\end{proof}

Thanks to Lemma 7.2 together with [V1, 3.1.Lemma], 
we see that the reduced von Neumann algebra $pMp$ 
is generated by
$$ 
\begin{aligned}
e^A_{1 i}\cdot(a\otimes 1)\cdot e_{j 1}^A &= 
\delta_{i j}\cdot(a\otimes p), \quad a \in A; \\
e^A_{1 i}\cdot e_{k \ell}^A\cdot e^A_{j 1} &= 
\delta_{i k}\cdot\delta_{\ell j}\cdot p; \\ 
e^A_{1 i}\cdot[u^n]_{k 1}^B\cdot e_{j 1}^A &= 
\delta_{i k}\cdot\delta_{j 1}\cdot u(n,i), \quad  
n \in {\mathbb Z},\ i = 2,3,\dots.  
\end{aligned}
$$
We set 
$$
N(n,i) := \{D_0\otimes{\mathbf C}p, u(n,i)\}'' \cong 
D_0\rtimes_{\alpha^n}{\mathbb Z} \quad 
\text{(thanks to Lemma 7.3)}
$$
with the conditional expectation 
$$
E_{(n,i)} = E_D^M|_{N(n,i)} : N(n,i) \rightarrow  
D_0\otimes{\mathbf C}p \cong D_0, 
$$
which coincides with the canonical one from 
$D_0\rtimes_{\alpha^n}{\mathbb Z}$ onto $D_0$. 

Summing up the discussions above, we conclude  

\begin{thm}
We have  
\begin{equation}
(pMp, (E_D^M)_p) \cong 
(A_0, E^{A_0}_{D_0}) *_{D_0} 
\left( \underset{\scriptstyle n \in {\mathbb Z} \atop 
\scriptstyle i = 2,3,\dots}{*}  
(N(n,i), E_{(n,i)})\right). 
\tag{7.6} 
\end{equation}
Here, the amalgamated free product 
$$
\underset{\scriptstyle n \in {\mathbb Z} \atop 
\scriptstyle i = 2,3,\dots}{*}  
(N(n,i), E_{(n,i)})
$$
is noting less than the crossed product of $D_0$ by the 
free group 
${\mathbb F}_{\infty}$ with countably many generators, 
whose action 
is defined as follows: 
\begin{equation}
(\text{\rm Id})^{*{\mathbb N}} * (\alpha)^{*{\mathbb N}} * 
(\alpha^2)^{*{\mathbb N}} * \cdots * (\alpha^n)^{*{\mathbb 
N}} * \cdots. 
\tag{7.7} 
\end{equation}
Here, $(\beta)^{*{\mathbb N}}$ means the free product of 
countably infinite 
copies of an automorphism
$\beta$. 
\end{thm}

\medskip
We further suppose that $A = B$, that is, $A_0 = B_0 = 
D_0 \rtimes_{\alpha}{\mathbb Z}$. Theorem 7.5 says 
that 
the reduced von Neumann algebra $pMp$ is isomorphic to 
the crossed product of $D_0$ by the free group ${\mathbb 
F}_{\infty}$ 
whose action is 
$$
\begin{aligned}
\alpha * &\left( (\text{\rm Id})^{*{\mathbb N}} * 
(\alpha)^{*{\mathbb N}} * (\alpha^2)^{*{\mathbb N}} * 
\cdots * (\alpha^n)^{*{\mathbb N}} * \cdots \right) \\ 
&\phantom{aaaaaaaaaa} = 
(\text{\rm Id})^{*{\mathbb N}} * (\alpha)^{*{\mathbb N}} * 
(\alpha^2)^{*{\mathbb N}} * \cdots * (\alpha^n)^{*{\mathbb 
N}} * \cdots. 
\end{aligned}
$$
Here, this equality follows from the simple fact: 
$\alpha * (\alpha)^{*{\mathbb N}} = (\alpha)^{*{\mathbb N}}$. 

\medskip\noindent
\begin{remark}  
{\rm The result obtained in this section is thought of 
as a negative evidence towards generalizing the work 
[G] on the invariant ``cost" of D. Gaboriau 
to general non-singular discrete measured groupoids. 
Roughly speaking, the ``cost" counts the number of 
free components in a given finite-measure preserving 
countable equivalence relation, and recently 
D. Shlyakhtenko [S2] generalized further to 
finite-measure preserving discrete groupoids from 
the free entropic viewpoint. 
Our result here says that the number of free components 
cannot be determined in the general non-singular case. 
Indeed, we suppose that our $A = B$ is a factor of type III 
(or of type II$_{\infty}$) and that $D$ is a Cartan 
subalgebra as before. Then the amalgamated free 
product $M$ is also a factor of type III and captured as 
a groupoid von Neumann algebra (see [Ks]). We can then 
choose an isometry $v \in A$ in such a way that 
$v D v^* = Dp$ with $vv^* = p \in D$. The $\text{Ad}v$ gives 
rise to an isomorphism between $M \supseteq D$ and 
$pMp \supseteq Dp$. Moreover, we can show 
\begin{equation}
(E_D^M)_p\circ\text{Ad}v = \text{Ad}v\circ E_D^M,  
\tag{7.8} 
\end{equation}
and hence $(M \supseteq D, E_D^M)$ can be identified 
with $(pMp \supseteq Dp, (E_D^M)_p)$, and the former 
has {\bf two} free components, but the latter has 
{\bf infinite} ones.} 
\end{remark}

\medskip
The discussions here (with trivial changes) also implies 

\begin{cor}
Let $N$ be an infinite injective 
factor of non-type {\rm I} with a Cartan subalgebra $D$. 
Then we have
\begin{equation}
N *_D N \cong N *_D N *_D N \cong \cdots 
\cong N *_D N *_D N *_D \cdots. 
\tag{7.9}  
\end{equation}
\end{cor}
 
\medskip
\begin{remark} {\rm One may replace $B(\ell^2({\mathbb N}))$ 
by $k\times k$ matrix algebra $M_k({\mathbf C})$ 
in the setting, and the discussion here still 
works without any essential change and the assertion 
(7.6) should be changed to 
\begin{equation}
(pMp, (E_D^M)_p) \cong 
(A_0, E^{A_0}_{D_0}) *_{D_0} 
\left( \underset{{\scriptstyle n \in {\mathbb Z} \atop 
\scriptstyle i = 2,3,\dots,k}}{*}  
(N(n,i), E_{(n,i)})\right) 
\tag{7.10}
\end{equation}
so that if $A_0 = B_0$ then $pMp$ is the crossed-product 
of $D_0$ by the free group ${\mathbb F}_{\infty}$ whose action 
is:
\begin{equation}
(\text{Id})^{*(k-1)} * (\alpha)^{*k} * 
(\alpha^2)^{*(k-1)} * \cdots * 
(\alpha^n)^{*(k-1)} * \cdots. 
\tag{7.11}
\end{equation}} 
\end{remark}

This is in particular thought of as 
a reduction formula of the amalgamated free product 
$R *_D R$ of two copies of the injective II$_1$ factor $R$ 
over a common Cartan subalgebra $D$ thanks to [CFW]. 
The result says that, if the amalgamated free product 
$R *_D R$ had the whole fundamental group (defined as in [P1]) 
${\mathcal F}(R *_D R \supseteq D) = {\mathbf R}_+^{\times}$, 
then the number of its free components would not be able to 
be determined uniquely. This is completely analogous to 
the situation of free group factors $L({\mathbb F}_n)$ with 
finite $n$ (see [V2, 6.13 {\it Remark}]). This analogy is 
very natural in a certain sense, because the amalgamated 
free product $R *_D R$ can be regarded as one candidate of 
the true generalizations of the free group factor 
$L({\mathbb F}_2)$ from the view-point of the idea 
generalizing the group von Neumann algebra construction 
to the group-measure space construction. (D. Shlyakhtenko 
[S1] provided another candidate ``$A$-valued semicircular 
systems" from the view-point of Voiculecu's 
free Gaussian functor (see [VDN]).) 

\bigskip

\section{References}

\medskip\noindent
[A] S. Adams, 
Boundary amenability for word hyperbolic groups and 
an application to smooth dynamics of simple groups. 
Topology {\bf 33} (1994), no. 4, 765--783. 

\medskip\noindent
[BD] E.F. Blanchard \& K.J. Dykema, 
Embeddings of reduced free products of operator algebras. 
preprint (1999). 

\medskip\noindent
[C1] A. Connes, 
Une classification des facteurs de type III.  
Ann. Scient. \`{E}c. Norm. Sup. {\bf 8} (1973), 133--252. 

\medskip\noindent 
[C2] A. Connes, 
Classification of injective factors. Cases II$_1$, II$_{\infty}$, 
III$_{\lambda}$, $\lambda \neq 1$. 
Ann. of Math. {\bf 104} (1976), 73--115. 

\medskip\noindent
[CFW] A. Connes, J. Feldman \& B. Weiss, 
An amenable equivalence relation is generated by a single transformation.
Ergod. Th. \& Dynam. Sys. {\bf 1} No. 4 (1982), 431--450.

\medskip\noindent
[CT] A. Connes \& M. Takesaki, 
The flow of weights on factor of type III.  
T\^{o}hoku Math. Journ., {\bf 29} (1977), 473--575. 

\medskip\noindent
[EFW] M. Enomoto, M. Fujii \& Y. Watatani, 
KMS states for gauge action on ${\mathcal O}_A$. 
Math. Japon. {\bf 29} (1984), no. 4, 607--619.
\medskip\noindent

\medskip\noindent
[FM] J. Feldman \& C.C. Moore, 
Ergodic equivalence relations, cohomology, and von Neumann algebras. 
I, II. 
Trans. Amer. Math. Soc., {\bf 234} (1977), 289--324, 325--358. 

\medskip\noindent
[G] D. Gaboriau, 
Co\^{u}t des relations d'\'{e}quivalence et des groupes. 
Invent. Math. {\bf 139} (2000), no. 1, 41--98.

\medskip\noindent 
[HOO] T. Hamashi, Y. Oka \& M. Osikawa, 
Flows associated with ergodic non-singular transformation groups. 
Publ. RIMS, Kyoto Univ. {\bf 11} (1975), 31--50. 

\medskip\noindent
[Kr] W. Krieger, 
On ergodic flows and the isomorphism of factors. 
Math. Ann. {\bf 223} (1976), 19--70. 

\medskip\noindent
[Ks] H. Kosaki, 
Free products of measured equivalence relations, 
Preprint. (2001) 

\medskip\noindent
[KS] G. Kuhn \& T. Steger, 
More irreducible boundary representations of free groups. 
Duke Math. J. {\bf 82} (1996), no. 2, 381--436.

\medskip\noindent
[P1] S. Popa, 
Some rigidity results in type II$_ 1$ factors. C. R. Acad. Sci. 
Paris S\'{e}r. I, Math. {\bf 311} (1990), no. 9, 535--538.

\medskip\noindent
[P2] S. Popa, 
Markov traces on universal Jones algebras and subfactors of 
finite index. 
Invent. Math. {\bf 111} (1993), 375--405. 

\medskip\noindent
[PS] C. Pensavalle \& T. Steger, 
Tensor products with anisotropic principal series representations 
of free groups. 
Pacific J. Math. {\bf 173} (1996), no. 1, 181--202.

\medskip\noindent
[RR] J. Ramagge \& G. Robertson, 
Factors from trees. 
Proc. Amer. Math. Soc. {\bf 125} (1997), no. 7, 2051--2055.

\medskip\noindent
[S1] D. Shlyakhtenko, 
$A$-valued semicercular systems. 
J. Funct. Annal. {\bf 166} (1999), 1--47. 

\medskip\noindent
[S2] D. Shlyakhtenko, 
Microstates free entropy and cost of equivalence relations. 
preprint (1999). 

\medskip\noindent
[Sp] J. Spielberg, 
Free-product groups, Cuntz-Krieger algebras, and covariant maps. 
Intern. J. Math. {\bf 2} (1991), 457--476. 

\medskip\noindent
[T1] M. Takesaki, 
Conditional expectations in von Neumann algebras. 
J. Funct. Annal. {\bf 9} (1972), 306--321. 

\medskip\noindent
[T2], M. Takesaki, 
Duality for crossed product and the structure of von Nenmann 
algebras of type III. 
Acta Math. {\bf 131} (1973), 249--310. 

\medskip\noindent
[U1] Y. Ueda, 
Amalgamated free product over Cartan subalgebra.  
Pacific J. Math. {\bf 191}, No.2 (1999), 359--392. 

\medskip\noindent
[U2] Y. Ueda,  
Fullness, Connes' $\chi$-groups, and ultra-products of 
amalgamated free products over Cartan subalgebras. 
Trans. Amer. Math. Soc. {\bf 355} (2003), 349-371 . 

\medskip\noindent
[Ver] A. M. Vershik, Trajectory Theory, chap. 5 in Dynamical 
Systems II, Ya. G. Sinai, ed. (Trandlated from Russian) 
(ENS) Encyclopaedia of Mathematical Sciences (Springer-Verlag, 
Berlin) {\bf 2}, 77-92 (1989)  

\medskip\noindent
[V1] D. Voiculescu, 
Circular and semicircular systems and free product factors. 
Operator algebras, unitary representations, enveloping algebras, 
and invariant theory (Paris, 1989), 45--60, Progr. Math. {\bf 92},
Birkh\"{a}user Boston, Boston, MA, 1990.

\medskip\noindent
[V2] D. Voiculescu, 
The analogues of entropy and Fisher's information measure in 
free probability theory, III: The absence of Cartan subalgebras. 
Geometric and Functional Analysis {\bf 6} (1996), 172--199. 

\medskip\noindent
[VDN] D.-V. Voiculescu, K.-J. Dykema \& A. Nica, 
Free Random Variables. 
CRM Monograph Series {\bf I}, Amer. Math. Soc. Providence, RI, 
1992. 

\medskip\noindent
[Z] R.J. Zimmer, 
Ergodic theory and semisimple groups.
Monographs in Mathematics, {\bf 81}. 
Birkh\"{a}user Verlag, Basel-Boston, Mass.,
1984. 

\enddocument